\documentclass{article}%
\usepackage{amsmath}
\usepackage{amsfonts}
\usepackage{amssymb}
\usepackage{graphicx}%
\providecommand{\U}[1]{\protect\rule{.1in}{.1in}}
\newtheorem{theorem}{Theorem}
\newtheorem{corollary}{Corollary}
\newtheorem{lemma}{Lemma}
\newtheorem{proposition}{Proposition}
\newenvironment{proof}[1][Proof]{\noindent\textbf{#1.} }{\ \rule{0.5em}{0.5em}}
\begin{document}

\title{Powers in finite groups}
\author{Nikolay Nikolov and Dan Segal}
\date{For Fritz Grunewald on his 60th birthday }
\maketitle

\begin{abstract}
If $G$ is a finitely generated profinite group then the verbal subgroup
$G^{q}$ is open. In a $d$-generator finite group every product of $q$th powers
is a product of $f(d,q)$ $q$th powers. 20E20, 20F20.

\end{abstract}

\section{Introduction}

\subsection{The main result}

For a group $H$ and positive integer $q$ the $q$th power subgroup is%
\[
H^{q}=\left\langle h^{q}\mid h\in H\right\rangle .
\]
Every element of $H^{q}$ is a product of $q$th powers; let us say that $H^{q}$
has \emph{width} $n$ if each such element is equal to a product of $n$ $q$th
powers (we don't assume that $n$ is minimal).

\begin{theorem}
\label{main}Let $q,~d\in\mathbb{N}$. Then there exists $f=f(d,q)$ such that
$H^{q}$ has width $f$ whenever $H$ is a $d$-generator finite group.
\end{theorem}

Straightforward arguments show that this is equivalent to

\begin{corollary}
\label{clos}If $G$ is a finitely generated profinite group and $q\in
\mathbb{N}$ then the (algebraically defined) subgroup $G^{q}$ has finite
width, and is closed in $G$.
\end{corollary}

Together with the positive solution of the Restricted Burnside Problem
(\cite{Z}, \cite{Z1}) this in turn implies

\begin{corollary}
\label{open}If $G$ is a finitely generated profinite group then $G^{q}$ is
open in $G$ for every $q\in\mathbb{N}$.
\end{corollary}

The deduction of the corollaries from Theorem \ref{main} is explained in
\cite{NS}, \S 1 and in Chapter 4 of \cite{W}. Theorem \ref{main} strengthens
\cite{NS}, Theorem 1.8 and Corollary \ref{open} generalizes \cite{NS}, Theorem 1.5.

For $q,~d\in\mathbb{N}$ let%
\[
\beta(d,q)
\]
denote the order of the $d$-generator restricted Burnside group of exponent
$q$; this is the maximal order of any finite $d$-generator of exponent
dividing $q$. The minimal size of a generating set for a group $H$ is denoted
$\mathrm{d}(H)$. If $H$ is finite and $\mathrm{d}(H)\leq d$ then $\left\vert
H:H^{q}\right\vert \leq\beta(d,q)$, so by Schreier's formula we have
$\mathrm{d}(H^{q})\leq d\beta(d,q)$. Taking%
\[
\delta(d,q)=d\beta(d,q)\cdot f(d,q)
\]
we see that Theorem \ref{main} implies

\begin{theorem}
\label{gen}Let $q,~d\in\mathbb{N}$. Then there exists $\delta=\delta(d,q)$
such that $H^{q}$ can be generated by $\delta$ $q$th powers in $H$ whenever
$H$ is a $d$-generator finite group.
\end{theorem}

\subsection{Wider implications}

The main results of \cite{NS} show that certain verbal subgroups are
necessarily closed in a finitely generated profinite group, namely those
associated to a locally finite word or to a simple commutator. This list can
now be extended:

\begin{theorem}
If $G$ is a finitely generated profinite group and $w$ is a non-commutator
word then the verbal subgroup $w(G)$ is open in $G$.
\end{theorem}

\noindent This greatly generalizes \cite{NS}, Theorem 1.3. It follows
immediately from Corollary \ref{open} since $w(G)$ contains $G^{q}$ where
$q=\left\vert \mathbb{Z}/w(\mathbb{Z})\right\vert $. Taking $G$ to be the free
profinite group on $d$ generators and $w$ any non-commutator word, we may
infer the existence of $f(d,w)$ and $\delta(d,w)$ such that if $H$ is any
$d$-generator finite group, then

\begin{itemize}
\item every product of $w$-values or their inverses in $H$ is equal to such a
product of length $f(d,w)$,

\item the verbal subgroup $w(H)$ is generated by $\delta(d,w)$ $w$-values
\end{itemize}

\noindent(cf. \cite{NS}, \S 1 or \cite{W}, \S 4.1).

Let us say that a group word $w$ is \emph{good} if $w(G)$ is closed in $G$
whenever $G$ is a finitely generated profinite group. The word $w=w(x_{1}%
,\ldots,x_{k})$ may be considered as an element of the free group $F$ on
$\{x_{1},\ldots,x_{k}\}$. Recall that $w$ is a \emph{commutator word} if $w\in
F^{\prime}$, the derived group of $F$. It is shown in \cite{JZ} that if $1\neq
w\in F^{\prime\prime}(F^{\prime})^{p}$ then $w(G)$ is \emph{not} closed in the
free pro-$p$ group $G$ on two generators ($p$ being any prime). Thus for a
non-trivial word $w$,%
\[
w\notin F^{\prime}~\Longrightarrow~w\text{ good}~\Longrightarrow~w\notin
F^{\prime\prime}(F^{\prime})^{p}~\forall p.
\]
The first implication is certainly strict, since simple commutators are good;
whether the second implication is reversible is an intriguing open question,
discussed at length in \cite{W}, Chapter 4.

\bigskip

This paper should be seen as a sequel to \cite{NS}, which contains all the
difficult arguments needed for Theorem \ref{main}. In particular, that paper
establishes (1) a weaker version of this theorem, restated below as
Proposition \ref{P1}, and (2) an implicit proof that Theorem \ref{main} would
follow from Theorem \ref{gen}; this is sketched in \S \ref{Sec3} below. As we
shall see, Theorem \ref{gen} can in turn be deduced quite easily from (1) and
another result in \cite{NS}.

The original motivation for \cite{NS} was to establish that every subgroup of
finite index in a finitely generated profinite group is open (`Serre's
problem'). This of course follows at once from Corollary \ref{open}, and our
initial strategy was indeed an attempt to prove the latter. Our failure to do
so forced us to develop machinery for dealing with other verbal subgroups;
this did the job just as well, and in fact better, in the sense that the
resulting proof was independent of the solution of the Restricted Burnside
Problem. Moreover, as far as we know, all the machinery of \cite{NS} is needed
to complete the proof of Theorem \ref{main}.

The main results all depend on the classification of finite simple groups,
which underpins much of \cite{NS}. The proof of Theorem \ref{main} also relies
on the solution of the Restricted Burnside Problem. This is inevitable:
indeed, Jaikin shows in \S 5.1 of \cite{JZ} how a positive solution to the
Restricted Burnside Problem for a prime power exponent $p^{n}$ can be deduced
directly from Corollary \ref{clos} with $q=p^{n+1}$.

Earlier special cases of Theorem \ref{main} were established in \cite{MZ},
\cite{SW} (for simple groups) and \cite{S} (for soluble groups).

\section{Preliminary results}

Henceforth all groups are assumed to be finite. We fix a positive integer $q$.
For a group $G$ and $m\in\mathbb{N}$ we write%
\begin{align*}
G_{q}  &  =\left\{  g^{q}\mid g\in G\right\} \\
G_{q}^{\ast m}  &  =\left\{  h_{1}h_{2}\cdots h_{m}\mid h_{1},\ldots,h_{m}\in
G_{q}\right\}  .
\end{align*}
Thus $G^{q}$ has width $m$ precisely when $G^{q}=G_{q}^{\ast m}$.

The largest integer $k$ such that $G$ involves the alternating group
$\mathrm{Alt}(k)$ as a section is denoted $\alpha(G)$.

\begin{proposition}
\label{P1}\emph{(\cite{NS}, Theorem 1.8)} Let $d,~k\in\mathbb{N}$. Then there
exists $h=h(k,d,q)$ such that $G^{q}$ has width $h$ whenever $G$ is a
$d$-generator finite group with $\alpha(G)\leq k$.
\end{proposition}

The next result is a slight weakening of \cite{NS}, Proposition 10.1:

\begin{proposition}
\label{P2}There exist $m=m(q)$ and $C(q)$ with the following property: if $N$
is a perfect normal subgroup of $G$ and $N/\mathrm{Z}(N)\cong S_{1}%
\times\cdots\times S_{n}$ where each $S_{i}$ is a non-abelian simple group
with $\left\vert S_{i}\right\vert >C(q)$ then%
\[
N\cdot G_{q}^{\ast m}=G_{q}^{\ast m}.
\]

\end{proposition}

We also need two simple lemmas. The first is a mild extension of a well-known
result due to Gasch\"{u}tz \cite{G} ; the proof given (for example) in
\cite{FJ}, Lemma 15.30 adapts easily to yield this version:

\begin{lemma}
\label{L1}Let $X\subseteq G$ and $N\vartriangleleft G$. Suppose that%
\[
G=N\left\langle X,y_{1},\ldots,y_{n}\right\rangle
\]
where $n\geq\mathrm{d}(G)$. Then there exist $a_{1},\ldots,a_{n}\in N$ such
that $G=\left\langle X,a_{1}y_{1},\ldots,a_{n}y_{n}\right\rangle $.
\end{lemma}

\begin{lemma}
\label{L2}Let $N\vartriangleleft G$. Then $G$ has a subgroup $L$ with $NL=G$
and $\alpha(L)\leq\max\{\alpha(G/N),4\}$.
\end{lemma}

\begin{proof}
Let $S$ be a Sylow $2$-subgroup of $N$ and put $L=\mathrm{N}_{G}(S)$. Then
$NL=G$ by the Frattini argument. If $\alpha(L)\geq5$ then $\alpha
(L)=\max\{\alpha(G/N),\alpha(L\cap N)\}$. The result follows since $L\cap N$
is an extension of a $2$-group by a group of odd order.
\end{proof}

\section{Generators}

Fix $k\geq5$ such that $k!>2C(q),$ and let $\mathcal{C}$ denote the class of
all groups $G$ with $\alpha(G)\leq k$. Put $m=m(q)$.

\begin{proposition}
\label{P3}Let $G$ be a $d$-generator group. Then $G=\left\langle X\cup
Y\right\rangle $ where $\left\vert X\right\vert \leq d$, $\left\vert
Y\right\vert \leq d$, $X\subseteq G_{q}^{\ast m}$ and $\left\langle
Y\right\rangle \in\mathcal{C}$.
\end{proposition}

\begin{proof}
Let $N$ be a minimal normal subgroup of $G$. Arguing by induction on the order
of $G$, we may suppose that $G=N\left\langle X^{\prime}\cup Y^{\prime
}\right\rangle $ where $\left\vert X^{\prime}\right\vert \leq d$, $\left\vert
Y^{\prime}\right\vert \leq d$, $X^{\prime}\subseteq G_{q}^{\ast m}$ and
$N\left\langle Y^{\prime}\right\rangle /N\in\mathcal{C}$. Applying Lemma
\ref{L2} to the group $N\left\langle Y^{\prime}\right\rangle $, we obtain a
set $Y^{\ast}$ with $\left\vert Y^{\ast}\right\vert =\left\vert Y^{\prime
}\right\vert $ such that $N\left\langle Y^{\ast}\right\rangle =N\left\langle
Y^{\prime}\right\rangle $ and $\left\langle Y^{\ast}\right\rangle
\in\mathcal{C}$. Then $G=N\left\langle X^{\prime}\cup Y^{\ast}\right\rangle $.
Say $X^{\prime}=\left\{  x_{1},\ldots,x_{d}\right\}  $ and $Y^{\ast}=\left\{
y_{1},\ldots,y_{d}\right\}  $ (allowing repeats if necessary).

\emph{Case 1.} Suppose that $N\notin\mathcal{C}$. By Lemma \ref{L1}, there
exist $a_{1},\ldots,a_{d}\in N$ such that $G=\left\langle Y^{\ast},a_{1}%
x_{1},\ldots,a_{d}x_{d}\right\rangle $. As $N\notin\mathcal{C}$, $N$ must be a
direct product of non-abelian simple groups of order exceeding $C(q)$. It
follows by Proposition \ref{P2} that $a_{i}x_{i}\in G_{q}^{\ast m}$ for each
$i$. The result follows with $X=\{a_{1}x_{1},\ldots,a_{d}x_{d}\}$, $Y=Y^{\ast
}$.

\emph{Case 2.} Suppose that $N\in\mathcal{C}$. Applying Lemma \ref{L1} again
we find $a_{1},\ldots,a_{d}\in N$ such that $G=\left\langle X^{\prime}%
,a_{1}y_{1},\ldots,a_{d}y_{d}\right\rangle $. Put $Y=\{a_{1}y_{1},\ldots
,a_{d}y_{d}\}$. Then $\left\langle Y\right\rangle \leq N\left\langle Y^{\ast
}\right\rangle \in\mathcal{C}$ and the result follows with $X=X^{\prime}$.
\end{proof}

\bigskip

We can now prove Theorem \ref{gen}. Let $H$ be a $d$-generator group.
According to Proposition \ref{P3},%
\[
H=\left\langle X\cup Y\right\rangle
\]
where $\left\vert X\right\vert \leq d$, $\left\vert Y\right\vert \leq d$,
$X\subseteq H_{q}^{\ast m}$ and $\left\langle Y\right\rangle \in\mathcal{C}$.
We apply Proposition \ref{P1} to the group $T=\left\langle Y\right\rangle $:
this shows that%
\[
T^{q}=T_{q}^{\ast h}%
\]
where $h=h(k,d,q)$. Put $\beta=\left\vert T:T^{q}\right\vert $; then
$\beta\leq\beta(d,q)$, and we have $T^{q}=\left\langle Z\right\rangle $ where
$\left\vert Z\right\vert \leq d\beta$.

Let $\{s_{1},s_{2},\ldots,s_{\beta}\}$ be a transversal to the cosets of
$T^{q}$ in $T$, put%
\begin{align*}
P  &  =\left\langle X\cup Z\right\rangle ,\\
K  &  =\left\langle P^{s_{1}},\ldots,P^{s_{\beta}}\right\rangle .
\end{align*}
Then $K\vartriangleleft H=KT$ and $\left\vert H:K\right\vert \leq\left\vert
T:T^{q}\right\vert =\beta$. Since $H^{q}=\left\langle H_{q}\right\rangle \geq
K$, it follows that $H^{q}=K\left\langle W\right\rangle $ for some subset $W$
of $H_{q}$ of size at most $\log_{2}\beta$.

Now each element of $Z$ is a product of $h$ $q$th powers in $T$ and each
element of $X$ is a product of $m$ $q$th powers in $H$; as $H^{q}$ is
generated by $W$ together with $\beta$ conjugates of $X\cup Z$, it follows
that $H^{q}$ can be generated by%
\[
\log_{2}\beta+\beta(dm+d\beta h)
\]
$q$th powers in $H$.

\section{Products of powers\label{Sec3}}

In the terminology of \cite{W}, Theorem \ref{gen} says that the word $x^{q}$
is $d$\emph{-restricted }for every $d$. Given this, Theorem \ref{main} becomes
a special case of \cite{W}, Theorem 4.7.9. However it seems worthwhile to make
this note self-contained modulo the paper \cite{NS}, so in this section we
sketch the deduction of Theorem \ref{main}.

This is an application of the main technical result of \cite{NS}; to state it
we need

\medskip

\noindent\textbf{Definition.} Let $G$ be a finite group and $K$ a normal
subgroup. Then $K$ is \emph{acceptable} if

\begin{enumerate}
\item[(i)] $K=[K,G]$ and

\item[(ii)] whenever $Z<N\leq K$ are normal subgroups of $G$, the factor $N/Z$
is not of the form $S$ or $S\times S$ for a non-abelian simple group $S$.
\end{enumerate}

\medskip

The `Key Theorem' stated in \cite{NS}, \S 2 is

\begin{proposition}
\label{P4}Let $K$ be an acceptable normal subgroup of $G=\left\langle
g_{1},\ldots,g_{\delta}\right\rangle $. Then%
\[
K=\left(
{\displaystyle\prod\limits_{i=1}^{\delta}}
[K,g_{i}]\right)  ^{\ast f_{1}}\cdot K_{q}^{\ast f_{2}}%
\]
where $f_{1}$ and $f_{2}$ depend only on $q$ and $\delta$.
\end{proposition}

(For a subset $X$ of $K$ we write $X^{\ast f}$ for the set $\left\{
x_{1}x_{2}\cdots x_{f}\mid x_{1},\ldots,x_{f}\in X\right\}  $.)

\bigskip

Let $H$ be a $d$-generator group and set $G=H^{q}$. As before, we have
$\mathrm{d}(G)\leq d^{\prime}=d\beta(d,q)$. Now $G$ has a series of
characteristic subgroups%
\[
K_{1}\geq K_{3}\geq K_{4}\geq K_{5}%
\]
such that

\begin{itemize}
\item $K_{5}$ is acceptable in $G$

\item $K_{3}$ is perfect and $K_{4}/K_{5}=\mathrm{Z}(K_{3}/K_{5})$

\item $K_{3}/K_{4}$ is a direct product of non-abelian simple groups of order
exceeding $C(q)$

\item $K_{1}/K_{3}$ is soluble

\item $\left\vert G:K_{1}\right\vert \leq\gamma=\gamma(d^{\prime},q)$
\end{itemize}

\noindent where $\gamma(d^{\prime},q)$ depends only on $d^{\prime}$ and $q$.
The proof, which is quite straightforward (given the classification of finite
simple groups), appears in \cite{NS}, \S 2 (see \emph{Proof of Theorem 1.6)}.

According to Theorem \ref{gen} there exist $g_{1},\ldots,g_{\delta}\in H_{q}$
such that $G=\left\langle g_{1},\ldots,g_{\delta}\right\rangle $ where
$\delta=\delta(d,q)$. Then $[h,g_{i}]\in H_{q}^{\ast2}$ for any $h\in H$ and
each $i$, so applying Proposition \ref{P4} we deduce that%
\[
K_{5}\subseteq H_{q}^{\ast(2\delta f_{1}+f_{2})}.
\]
Proposition \ref{P2} shows that $K_{3}\subseteq H_{q}^{\ast m}\cdot K_{5}$.
Now let $k^{\prime}\geq\max\{5,q+2\}$ be such that $k^{\prime}!>2\gamma
(d^{\prime},q)$. Then $\alpha(H/K_{3})\leq k^{\prime}$; thus Proposition
\ref{P1} gives%
\[
H^{q}\subseteq H_{q}^{\ast h}\cdot K_{3}%
\]
where $h=h(k^{\prime},d,q)$. Putting everything together we get $H^{q}%
\subseteq H_{q}^{\ast f}$ where%
\[
f=h+m+2\delta f_{1}+f_{2},
\]
a number that depends only on $d$ and $q$. This completes the proof of Theorem
\ref{main}.

\bigskip

N. Nikolov

Dept. of Mathematics

Imperial College

London SW7 2AZ

\texttt{n.nikolov@imperial.ac.uk}

\medskip

D. Segal

All Souls College

Oxford OX1 4AL

\texttt{dan.segal@all-souls.ox.ac.uk}

\end{document}